\documentclass[invmat,referee]{svjour} 
\usepackage{amssymb,amsmath,latexsym,epsfig,euscript} 

\newlength{\jmr}
\newlength{\hwl}
\newlength{\khov}
\newlength{\bernd}
\newtheorem{smi}{Smirnov's Theorem}

\newtheorem{thm}{Theorem}

\newtheorem{dfn}{Definition}
\newtheorem{cor}{Corollary}
\newtheorem{rem}{Remark}	
\newtheorem{rems}[rem]{Remarks} 
\newtheorem{ex}{Example}

\renewcommand{\mod}{\mathbf{mod}}

\newcommand{\np}{{\mathbf{NP}}}

\newcommand{\bpp}{{\mathbf{BPP}}}

\newcommand{\pp}{\mathbf{P}}

\newcommand{\eps}{\varepsilon}

\newcommand{\cO}{\mathcal{O}}
\newcommand{\supp}{\mathrm{Supp}}

\newcommand{\conv}{\mathrm{Conv}}

\newcommand{\thth}{{\underline{\mathrm{th}}}}

\newcommand{\nd}{{\underline{\mathrm{nd}}}}

\newcommand{\F}{\mathbb{F}}
\newcommand{\Q}{\mathbb{Q}}
\newcommand{\R}{\mathbb{R}}
\newcommand{\C}{\mathbb{C}}
\newcommand{\N}{\mathbb{N}}
\newcommand{\Z}{\mathbb{Z}}
 
\newcommand{\cp}{\mathfrak{p}}

\newcommand{\ord}{\mathrm{ord}}

\newcommand{\newt}{\mathrm{Newt}}

\newcommand{\Zn}{\Z^n}

\newcommand{\Rn}{\R^n}

\newcommand{\Cn}{\C^n}

\newcommand{\Csn}{{(\C^*)}^n}
\renewcommand{\qed}{$\blacksquare$}
\newcommand{\dia}{$\diamond$}
\newcommand{\cF}{{\mathcal{F}}}
\newcommand{\cL}{{\mathcal{L}}}
\newcommand{\cM}{{\mathcal{M}}}

\newcommand{\bO}{\mathbf{O}}

\begin{document}

\begin{center}
{\bf ARITHMETIC MULTIVARIATE DESCARTES' RULE}

{\bf J. Maurice Rojas} 

Department of Mathematics 

Texas A\&M University 

TAMU 3368

College Station, TX \ 77845-3368

USA

\end{center}

\title{Arithmetic Multivariate Descartes' Rule} 

% \subjclass{Primary 
% 11G25; % Arithmetic Algebraic Geometry; varieties over finite and 
%        % local ground fields
% Secondary 
% 11G35, % ...global ground fields
% 14D10, % algebraic geometry; families, fibrations; arithmetic ground fields
% 14G20. % algebraic geometry; arithmetic problems/diophantine geometry; 
%        % local ground fields 
% } 

\author{J.\ Maurice Rojas}

\thanks{
This research was partially supported by a grant from the 
Texas A\&M College of Science.} 

\institute{Department of Mathematics, Texas A\&M University, TAMU 3368, 
College Station, Texas 77843-3368, USA. }  
\email{rojas@math.tamu.edu \hfill {\rm Web Page:} 
http://www.math.tamu.edu/\~{}rojas\\
\mbox{}\hfill {\rm FAX:} (979) 845-6028 \hfill\mbox{}} 

\date{November 15, 2001} 

\maketitle

\begin{abstract} 
Let $\cL$ be any number field or $\cp$-adic field and consider 
$F\!:=\!(f_1,\ldots,f_k)$ where $f_i\!\in\!\cL[x^{\pm 1}_1,\ldots,
x^{\pm 1}_n]\setminus\{0\}$ for all $i$ 
and there are exactly $m$ distinct exponent vectors appearing in 
$f_1,\ldots,f_k$. 
We prove that $F$ has no more than $1+\left(\sigma mn(m-1)^2\log m\right)^n$ 
geometrically isolated roots in $\cL^n$, where $\sigma$ is an explicit 
and effectively computable constant depending only on $\cL$. This gives a 
significantly sharper arithmetic analogue of 
Khovanski's Theorem on Fewnomials and a higher-dimensional generalization of 
an earlier result of Hendrik W.\ Lenstra, Jr.\ for the case of a single 
univariate polynomial. We also present some further refinements of our new 
bounds and briefly \mbox{discuss the complexity of finding isolated rational 
roots.}  
\end{abstract} 

\section{Introduction} 
\label{sec:intro} 
A consequence of Descartes' Rule (a classic result dating back to 1637) 
is that any real univariate polynomial with exactly $m\!\geq\!1$ monomial 
terms has at most $2m-1$ real roots. This has since been generalized 
by Askold G.\ Khovanski during 1979--1987 (see \cite{kho} and 
\cite[pg.\ 123]{few}) to certain systems of multivariate sparse
polynomials and even {\bf fewnomials}.\footnote{Sparse 
polynomials are sometimes also known as lacunary polynomials and, over $\R$, 
are a special case of fewnomials --- a more general class of real 
analytic functions \cite{few}. } 
Here we provide ultrametric and thereby arithmetic analogues for both 
results: we give explicit upper bounds, independent of the degrees of the 
underlying polynomials, for the number of isolated roots of sparse 
polynomial systems over any {\bf $\pmb{\cp}$-adic field} and, 
as a consequence, over any {\bf number field}. For convenience, let us 
henceforth respectively refer to these cases as the {\bf local} 
case and the {\bf global} case. 

Suppose $\pmb{f_1},\ldots,\pmb{f_k}\!\in\!\cL[x^{\pm 1}_1,\ldots,
x^{\pm 1}_n]\setminus\{0\}$ 
where $\pmb{\cL}$ is a field to be specified later, and
$\pmb{m}$ is the total number of distinct exponent vectors appearing 
in $f_1,\ldots,f_k$ (assuming all polynomials are written as sums of 
monomials). We call $\pmb{F}\!:=\!(f_1,\ldots,f_k)$ an {\bf $\pmb{m}$-sparse 
$\pmb{k\times n}$ polynomial system over $\pmb{\cL}$}. 
Khovanski's results take $\cL\!=\!\R$ 
and yield an explicit upper bound for the number of non-degenerate roots, in 
the non-negative orthant, of any $m$-sparse $n\times n$ polynomial system 
\cite{kho,few}. With a little extra work (e.g., \cite[cor.\ 3.2]{real}) his 
results imply an upper bound of $2^{\cO(n)}n^{\cO(m)}2^{\cO(m^2)}$ on the 
number of isolated\footnote{We say a root of $F$ is {\bf geometrically 
isolated} iff it is a zero-dimensional component of the underlying 
scheme over the algebraic closure of $\cL$ defined by $F$. For the case 
of $\cL\!=\!\R$ one can in fact use the slightly looser definition 
that a point is {\bf topologically} isolated iff it is a connected 
component of the underlying 
real zero set. Unless otherwise mentioned, all our isolated roots will 
be geometrically isolated.} roots of $F$ in $\Rn$, and this is asymptotically 
the best general upper bound currently known. In particular, since 
it is easy to show that the last bound can in fact be replaced by 
$1$ when $m\!\leq\!n$ (see, e.g., \cite[thm.\ 3, part (b)]{tri}), 
one should focus on better understanding the behavior of the 
maximum number of isolated real 
roots for $n$ fixed and $m\!\geq\!n+1$. For example, is the dependence on $m$  
in fact polynomial for fixed $n$? This turns out to be an open 
question, but we can answer the arithmetic analogue (i.e., where 
$\cL$ is any $\cp$-adic field or any number field) affirmatively and 
explicitly: 
\begin{thm}
\label{thm:big}
Let $p$ be any (rational) prime and 
$d,\delta$ positive integers. Suppose 
$\cL$ is any degree $d$ algebraic extension of $\Q_p$ or $\Q$, and 
let $\cL^*\!:=\!\cL\setminus\{0\}$. Also let $F$ 
be any $m$-sparse $k\!\times\!n$ polynomial system over $\cL$ and 
define $\pmb{B(\cL,m,n)}$ to be 
the maximum number of isolated roots in $(\cL^*)^n$ of such an 
$F$ in the local case, counting multiplicities\footnote{ The 
multiplicity of any isolated root here, which we take in the sense of 
intersection theory for a scheme over the algebraic closure of $\cL$ 
\cite{fulton}, turns out to always be a positive integer (see, e.g., 
\cite{smirnov,jpaa}).} if (and only if) $k\!=\!n$. 
Then $B(\cL,m,n)\!=\!0$ (if $m\!\leq\!n$ or 
$k\!<\!n$) and\\ 
\scalebox{.97}[1]{$B(\cL,m,n)\!\leq\!u(m,n)\left\{c(m-1)n(p^d-1)\left[1
+d\log_p\left(\frac{d(m-1)}{\log p}\right)\right]\right\}^n$ 
(if $m\!\geq\!n+1$ and $k\!\geq\!n$),}\\ 
where $u(m,n)$ is $m-1$, $4(m-1)^2$, or $(m(m-1)/2)^n$ according as 
$n\!=\!1$, $n\!=\!2$, or $n\!\geq\!3$; 
$c\!:=\!\frac{e}{e-1}\!\leq\!1.582$ and 
$\log_p(\cdot)$ denotes the base $p$ logarithm function. 
Furthermore, moving to 
the global case, let us say a root $x\!\in\!\C^n$ of $F$ is of {\bf degree 
$\pmb{\leq\!\delta}$ over $\cL$} iff every coordinate of $x$ lies in an 
extension of degree $\leq\!\delta$ of $\cL$, and let us define 
$\pmb{A(\cL,\delta,m,n)}$ to be
the maximum number of isolated roots of $F$ in $\Csn$ of degree 
$\leq\!\delta$ over $\cL$, counting multiplicities$^3$ if (and only if) 
$k\!=\!n$. Then 
$A(\cL,\delta,m,n)\!=\!0$ (if $m\!\leq\!n$ or $k\!<\!n$) and\\ 
\scalebox{.95}[1]{$A(\cL,\delta,m,n)\!\leq\!2u(m,n)\left\{c(m-1)n 2^{d\delta}
\left[1+2d^2\delta^2\log_2\left(\frac{d^2\delta^2(m-1)}{\log 2}\right)
\right]\right\}^n$ (if $m\!\geq\!n+1$ and $k\!\leq\!n$).} 
\end{thm} 

Our bounds can be sharpened even further, and this is detailed 
in corollaries \ref{cor:local} and \ref{cor:global} of sections \ref{sec:local} 
and \ref{sec:global}, respectively. 
\begin{rem} 
\label{rem:aff} 
At the expense of underestimating\footnote{e.g., 
roots on the coordinate hyperplanes 
may have multiplicities $>\!1$ counted as $1$ instead.} some multiplicities, 
we can easily obtain a bound for the number of isolated roots of $F$ in 
$\cL^n$ (in the local case) or roots in $\Cn$ of degree $\leq\!\delta$ over 
$\cL$ (in the global case): By simply setting all possible subsets of 
variables to zero, we easily obtain respective bounds of 
$1+\sum^n_{j=1}\begin{pmatrix}n\\ j 
\end{pmatrix}B(\cL,m,j)\!\leq\!1+2^nB(\cL,m,n)$ and 
$1+\sum^n_{j=1}\begin{pmatrix}n\\ j 
\end{pmatrix}A(\cL,\delta,m,j)\!\leq\!1+2^nA(\cL,\delta,m,n)$. 
Of course, since many of the monomial terms of 
$F$ will vanish upon setting an $x_i$ to $0$, these bounds will usually be 
larger than really necessary. \dia 
\end{rem}  

\begin{ex} 
\label{ex:3m} 
Consider the following $2\times 2$ system over $\Q_2$: 
\begin{center}
$f_1(x_1,x_2)\!:=\!\alpha_1+\alpha_2x^{u_1}_1x^{u_2}_2
+\alpha_3x^{v_1}_1x^{v_2}_2$

$f_2(x_1,x_2)\!:=\!\beta_1+\beta_2x^{a_{2,1}}_1x^{a_{2,2}}_2
+\cdots +\beta_\mu x^{a_{\mu,1}}_1x^{a_{\mu,2}}_2$
\end{center} 
which is $m$-sparse for some $m\!\leq\!\mu+2$. Theorem \ref{thm:big}  
and an elementary calculation then tell us that such an $F$ has no more 
than 
\begin{center} 
$41(\mu+1)^4\left(1+\log_2\left(\frac{\mu+1}{0.693}
\right)\right)^2$
\end{center} 
isolated roots, counting multiplicities, in $(\Q^*_2)^2$  
(and $(\Q^*)^2$ as well, via the natural embedding $\Q\hookrightarrow \Q_2$). 
For instance, $\mu\!=\!3 \Longrightarrow F$ is at worst $5$-sparse and 
has\footnote{The numerical calculations throughout this paper were 
done with the assistance of {\tt Maple}, and the code for these 
calculations is available from the author's web-page. }  
no more than $127645$ roots in $(\Q^*_2)^2$. Explicit bounds 
independent of the total degrees of $f_1$ and $f_2$ appear to have been 
unknown before. 
However, if we replace $\Q_2$ by $\R$ throughout, then the best previous
upper bounds were  $4(2^\mu-2)$ for all $\mu\!\geq\!1$ 
and a bound of $20$ in the special case $\mu\!=\!3$. Interestingly, the latter 
bounds, which follow easily from \cite[thm.\ 1]{tri}, in fact allow us to 
take real {\bf exponents} and count {\bf topologically} isolated roots, but 
without multiplicities.\footnote{
Khovanski's Theorem on Fewnomials \cite[cor.\ 7, sec.\ 3.12]{few},
which only counts roots with non-singular Jacobian, implies an upper bound of   $995328$ for $\mu=3$. } The real analytic upper bound exceeds our 
arithmetic bound for all $\mu\!\geq\!29$, where both bounds begin 
to exceed $1.3$ billion. A sharper bound, based on a 
refinement of theorem \ref{thm:big} (cf.\ corollary \ref{cor:local}),  
appears in example \ref{ex:3m2} of section \ref{sec:local}. \dia 
\end{ex}  
\begin{ex} 
Another consequence of theorem \ref{thm:big} is that for {\bf fixed} $\cL$, 
we now know $\log B(\cL,m,n)$ and 
$\log A(\cL,m,n)$ to within a constant factor: For $m\!\geq\!2$ consider the 
$m$-sparse $n\times n$ polynomial system $F\!=\!(f_1,\ldots,f_n)$ where 
$f_i\!=\!\prod^{m-1}_{j=1}\left(x_i-j\right)$ for all $i$. 
Clearly then, this $F$ has exactly $(m-1)^n$ 
isolated roots in $\N^n$. It is curious that the  
analogous growth-rate is unknown if $\cL$ is replaced 
by the usual Archimedean \mbox{completion $\R$ of $\Q$. \dia} 
\end{ex} 

A weaker version of theorem \ref{thm:big} with non-explicit bounds  
was derived earlier in \cite{finite}. In particular, explicit bounds 
were known previously only in the case $k\!=\!n\!=\!1$ \cite[thm.\ 1 and 
thm.\ 2]{lenstra2}, and all our bounds (save the global case) match the  
bounds of \cite{lenstra2} in this special case.\footnote{In order to 
streamline the proof of our number field generalization, we left our bound on 
$A(\cL,\delta,m,n)$ in theorem \ref{thm:big} a bit loose: for $n\!=\!1$ our 
bound reduces to $\cO\left(d^2\delta^2 m^2 2^{d\delta} \log(d\delta m)
\right)$, while the older univariate result yields 
$\cO\left(d\delta m^2 2^{d\delta} \log(d\delta m) \right)$ in our notation. 
A sharper bound, agreeing with Lenstra's univariate bound when $n\!=\!1$, 
appears in corollary \ref{cor:global} of section \ref{sec:global}. } 
Philosophically, the approach of \cite{lenstra2} was more algebraic (low 
degree factors of polynomials) while our point of view here is more geometric 
(isolated rational points of low degree in a hypersurface 
intersection).\footnote{
Lenstra has also considered a higher-dimensional generalization but in a 
different direction: bounds for the number of rational hyperplanes 
in a hypersurface defined by a single $m$-sparse $n$-variate polynomial 
\cite[prop.\ 6.1]{lenstra2}.} The only 
other results known for $k\!>\!1$ or $n\!>\!1$ were derived via rigid analytic 
geometry and model theory, and in our notation yield a non-effective bound of 
$B(\Q_p,m,n)\!<\!\infty$ (see the seminal works \cite{vandenef,lipshitz}). 

Our approach is simpler and is based on  
a higher-dimensional generalization (theorem \ref{thm:lip} of the next 
section) of a result of Hendrik W.\ Lenstra, Jr.\ for univariate sparse 
polynomials over certain algebraically closed fields \cite[thm.\ 3]{lenstra2}. 
Indeed, aside from the introduction of some higher-dimensional 
convex geometry, our proof of theorem \ref{thm:big} 
is structurally quite similar to Lenstra's proof of the $k\!=\!n\!=\!1$ case 
in \cite{lenstra2}: reduce the global case to the local case, then reduce 
the local case to a refined result over the $p$-adic {\bf complex} numbers. 

We now describe two results used in our proofs which may be of 
broader interest. We also point out that connections between 
our results and complexity theory, including the question of 
whether we can find isolated rational roots in polynomial time, is 
described in section \ref{sec:pnp}. 

\subsection{The Distribution of $p$-adic Complex Roots}  
\label{sec:id} 
For any (rational) prime $p$, let $\C_p$ denote the completion (with respect 
to the $p$-adic metric) of the algebraic closure of $\Q_p$. 
Theorem \ref{thm:big} follows from a careful application of 
two results on the distribution of roots of $F$ in $(\C^*_p)^n$. 
The first result strongly limits the number of roots that 
can be $p$-adically close to the point $(1,\ldots,1)$. The second result 
strongly limits the number of distinct valuation vectors which can occur 
for the roots of $F$. 
\begin{thm} 
\label{thm:lip} 
Let $F$ be any $m$-sparse $k\times n$ polynomial system over $\C_p$. 
Also let $r_1,\ldots,r_n\!>\!0$, $\pmb{r}\!:=\!(r_1,\ldots,r_n)$, 
and let $\ord_p : \C_p \longrightarrow \Q\cup\{+\infty\}$ denote the 
usual exponential valuation, normalized\footnote{So, for example, 
$\ord_p 0\!=\!+\infty$ and $\ord_p(p^kr)\!=\!k$ whenever $r$ is a unit in 
$\Z_p$ and $k\!\in\!\Q$.} so that $\ord_p p \!=\!1$. Finally, let 
$\pmb{C_p(m,n,r)}$ denote the maximum number of isolated roots 
$(x_1,\ldots,x_n)$ of $F$ in $\C^n_p$ with $\ord_p(x_i-1)\!\geq\!r_i$ 
for all $i$, counting multiplicities if (and only if) $k\!=\!n$. Then 
$C_p(m,n,r)\!=\!0$ (if $m\!\leq\!n$ or $k\!<\!n$) 
and \[C_p(m,n,r)\!\leq\!\left.\left\{c(m-1)\left[r_1+\cdots+ r_n +\log_p
\left(\frac{(m-1)^n}{r_1\cdots r_n\log^n p}\right)\right]\right\}^n
\right/\prod^n_{i=1}r_i\] (if $m\!\geq\!n+1$ and $k\!\geq\!n$), 
where $c\!:=\!\frac{e}{e-1}\!\leq\!1.582$. Furthermore, if we restrict to 
those $F$ where $k\!=\!n$ and $f_i$ has exactly $m_i$ monomial terms for all 
$i$, then we have the sharper bounds of $C_p(m,n,r)\!=\!0$ (if $m_i\!\leq\!1$ 
for some $i$) and 
\[C_p(m,n,r)\!\leq\!c^n\prod\limits^n_{i=1}\left\{\left. (m_i-1)
\left[r_1+\cdots+ r_n +\log_p\left(\frac{(m_i-1)^n}{r_1\cdots r_n\log^n p}
\right)\right]\right/r_i\right\}\] (if $m_1,\ldots,m_n\!\geq\!2$). 
\end{thm} 
\noindent 
These bounds appear to be new: the only previous results in this 
direction appear to have been Lenstra's derivation of 
the special case $n\!=\!1$ \cite[thm.\ 3]{lenstra2} and an earlier 
observation of Leonard Lipshitz \cite[thm.\ 2]{lipshitz} equivalent to the 
non-explicit bound \mbox{$C_p(m,n,(1,\ldots,1))\!<\!\infty$.} 

Our last bound over $\C^n_p$ is based on a toric arithmetic-geometric 
result of Smirnov, stated below. 
\begin{dfn} 
\label{dfn:poly}
For any $\pmb{a}\!=\!(a_1,\ldots,a_n)\!\in\!\Zn$, let 
$\pmb{x^a}\!:=\!x^{a_1}_1\cdots x^{a_n}_n$. 
Writing any $\pmb{f}\!\in\!\cL[x_1,\ldots,x_n]$ as $\sum_{a\in \Zn} 
c_ax^a$, we call $\pmb{\supp(f)}\!:=\!\{a \; | \; c_a\!\neq\!0\}$ the 
{\bf support} of $f$. Also, let $\pmb{\pi} : \R^{n+1}\longrightarrow \Rn$ be 
the natural projection forgetting the $x_{n+1}$ coordinate and, for 
any $n$-tuple of polytopes $P\!=\!(P_1,\ldots,P_n)$, 
define $\pmb{\pi(P)}\!:=\!(\pi(P_1),\ldots, \pi(P_n))$. \dia 
\end{dfn} 
\begin{dfn} 
For any $k\times n$ polynomial system $F$ over $\C_p$, its {\bf 
$\pmb{k}$-tuple of $\pmb{p}$-adic Newton polytopes}, 
$\pmb{\newt_p(F)}\!=\!(\newt_p(f_1),\ldots,\newt_p(f_k))$, is 
defined as follows: $\pmb{\newt_p(f_i)}\!:=\!\conv(\{(a,\ord_p c_a) \; | \; 
a\!\in\!\supp(f_i) \})\!\subset\!\R^{n+1}$, where $\pmb{\conv(S)}$ denotes the 
convex hull of\footnote{i.e., smallest convex set containing...} 
a set $S\!\subseteq\!\R^{n+1}$. 
Also, for any $w\!\in\!\Rn$ and any 
closed subset $B\!\subset\!\Rn$, let the {\bf face of $\pmb{B}$ with inner 
normal $\pmb{w}$}, $\pmb{B^w}$, be the set of points $x\!\in\!B$ which 
minimize the inner product $w\cdot x$. Finally, 
let $\newt^w_p(F)\!:=\!(\newt^w_p(f_1),\ldots,\newt^w_p(f_k))$. \dia 
\end{dfn} 

\begin{smi} 
\cite[thm.\ 3.4]{smirnov}
\label{thm:smirnov} 
For any $n\!\times\!n$ polynomial system $F$ over $\C_p$, 
the number of isolated roots $(x_1,\ldots,x_n)$ of $F$ in $(\C^*_p)^n$ 
satisfying $\ord_p x_i\!=\!r_i$ for all $i$ (counting multiplicities) 
is no more than $\cM(\pi(\newt^{\hat{r}}_p(F)))$, \mbox{where 
$\pmb{\hat{r}}\!:=\!(r_1,\ldots,r_n,1)$, 
$\pmb{\cM(\cdot)}$ denotes {\bf mixed}} {\bf volume} \cite{buza} 
(normalized so that \scalebox{.9}[1]{$\cM(\conv(\{\bO,e_1,\ldots,e_n\}),\ldots,
\conv(\{\bO,e_1,\ldots,e_n\}))\!=\!1$),} and 
$\pmb{e_i}$ is the $i^\thth$ standard basis vector of $\Rn$. \qed 
\end{smi} 

\begin{rems}
\label{rems:convex} 
\begin{itemize}
\item[]{}
\item[0.]{ Explicit examples of the preceding constructions are  
illustrated in \cite{finite}.}  
\item[1.]{ The number of possible distinct valuation vectors for a root of 
an $n$-variate polynomial system $F$ can thus be combinatorially bounded 
from above as a function depending solely on $n$ and the number of monomial 
terms (cf.\ section \ref{sec:local}).} 
\item[2.]{ The number of roots of $F$ in $(\C^*_p)^n$ with given 
valuation vector thus depends strongly on the individual exponents of $F$ --- 
not just on the number of monomial terms.}
\item[3.]{It is thus only the {\bf lower}\footnote{Those with positive 
$x_{n+1}$ coordinate for their inner normals...} faces of the $p$-adic 
Newton polytopes that matter in counting roots or valuation vectors. \dia} 
\end{itemize}
\end{rems}

We prove theorem \ref{thm:lip} in section \ref{sec:lip}. However, 
let us first show how theorem \ref{thm:lip} implies theorem 
\ref{thm:big}: We will begin by examining the local case in the 
next section, and then complete our proof by deriving the global 
case from the local case in section \ref{sec:global}. 

\section{The Local Case of Theorem \ref{thm:big}} 
\label{sec:local} 
Here we will assume that $\cL$ is any degree $d$ algebraic extension 
of $\Q_p$. The following lemma will help us reduce to the case 
$k\!=\!n$. 
\begin{lemma}
\label{lemma:gh}
(See \cite[lemma 1]{finite}.) 
Following the notation of theorem \ref{thm:big}, there is a matrix 
$[a_{ij}]\!\subset\!\Z^{n\times k}$ such that the zero set of 
$G\!:=\!(a_{11}f_1+\cdots+a_{1k}f_k,\ldots,a_{n1}f_1+\cdots+a_{nk}f_k)$
in $\Cn_p$ is the union of the zero set of $F$ in $\Cn_p$ and a finite
(possibly empty) set of points. \qed 
\end{lemma}  

\noindent
{\bf Proof of the Local Case of Theorem \ref{thm:big}:} It is 
clear that there are no isolated roots whatsoever if $k\!<\!n$, since 
the underlying algebraic set over $\C^n_p$ is positive-dimensional. 
So we can assume $k\!\geq\!n$. In the event that $k\!>\!n$, lemma 
\ref{lemma:gh} then allows us to replace $F$ by a new 
$n\times n$ polynomial system (with no new exponent vectors) which 
has at least as many isolated roots as our original $F$. 
So we can assume $k\!=\!n$ and observe that root multiplicities 
are preserved if lemma \ref{lemma:gh} was not used (i.e., 
if we already had $k\!=\!n$ in our input). Since an $m$-sparse $n\times n$ 
polynomial system clearly has no geometrically isolated roots whatsoever 
when $m\!\leq\!n$, we can clearly assume 
that $m\!\geq\!n+1$. (Indeed, upon dividing each equation by a suitable 
monomial, $m\!\leq\!n$ clearly implies that we can obtain $n$ linear equations 
in $\leq n-1$ non-constant monomial terms.)  

The well-known classification of when mixed volumes vanish 
\cite{buza} then yields that $\cM(\pi(\newt^{\hat{r}}_p(F)))\!>\!0$ 
$\Longrightarrow$ there are linearly independent vectors 
$v_1,\ldots,v_n$, with $v_i$ an edge of  
$\newt^{\hat{r}}_p(f_i)$ for all $i$. So let $\eps_i$ be the number 
of edges of $\newt_p(f_i)$. If $n\!=\!1$ then we clearly have 
$\eps_i\!\leq\!m-1$ for all $i$, and this is a sharp bound 
for all $m$. If $n\!>\!2$ then we have the 
obvious bound of $\eps_i\!\leq\!m(m-1)/2$ 
for all $i$, and it is not hard to generate examples showing that this 
bound is sharp for all $m$ as well \cite[thm.\ 6.5, pg.\ 101]{ede}. 
If $n\!=\!2$ then note that the number of edges of $\newt_p(f_i)$ is 
clearly not decreased if we triangulate the boundary of $\newt_p(f_i)$. 
Since each $2$-face of the resulting complex is incident to exactly 
$2$ edges, Euler's relation \cite[thm.\ 6.8, pg.\ 103]{ede} then 
immediately implies that $\eps_i\!\leq\!2m-2$ for all $i$. (This bound is 
easily seen to be sharp for all $m\!\geq\!4$.) 

We thus obtain that there are no more than $\eps_1\cdots \eps_n\!\leq\!u(m,n)$ 
possible values for an $r\!\in\!\Rn$ with $\hat{r}\!=\!(r,1)$ and 
$\cM(\pi(\newt^{\hat{r}}_p(F)))\!>\!0$.  
In particular, by Smirnov's Theorem, this implies that 
the number of distinct values for the valuation vector 
$(\ord_p x_1,\ldots,\ord_p x_n)$, where 
$(x_1,\ldots,x_n)\!\in\!(\C^*_p)^n$ 
is a root of $F$, is no more than $u(m,n)$. So let us temporarily fix 
$(\pmb{r_1},\ldots,\pmb{r_n})\!:=\!r$ and see how many roots of $F$ in 
$(\cL^*)^n$ can have valuation vector $r$. 

Following the notation of theorem \ref{thm:lip}, let $R_p\!:=\!\{x\!\in\!\C_p 
\; | \; |x|_p\!\leq\!1\}$ be the ring of algebraic integers 
of $\C_p$, let $M_p\!:=\! \{x\!\in\!\C_p \; | \; |x|_p\!<\!1\}$ 
be the unique maximal ideal of $R_p$, \mbox{$\F_\cL\!:=\!(R_p\cap \cL)
/(M_p\cap \cL)$,}  
and let $\rho$ be any generator of the principal ideal $M_p\cap \cL$ 
of $R_p\cap \cL$. Also let $e_\cL\!:=\!\max_{y\in \cL^*}\{ 
|\ord_p y|^{-1}\}$ and $q_\cL\!:=\!\#\F_\cL$. 
(The last two quantities are respectively known as the {\bf ramification 
degree} and {\bf residue field cardinality} of $\cL$, and satisfy 
$e_\cL,\log_p q_\cL\!\in\!\N$ and $e_\cL\log_p q_\cL=d$ 
\cite[ch.\ III]{koblitz}.) Since $\ord_p \rho\!=\!1/e_\cL$, it 
is clear that $r$ a valuation vector of a root of $F$ in 
$(\cL^*)^n \Longrightarrow r\!\in\!(\Z/e_\cL)^n$. 

Fixing a set $A_\cL\!\subset\!R_p$ of representatives for $\F_\cL$ 
(i.e., a set of $q_\cL$ elements of $R_p\cap \cL$, exactly {\bf one} of 
which lies in $M_p$, whose image 
mod $M_p\cap \cL$ is $\F_\cL$), we can then write any 
$x_i\!\in\!\cL$ uniquely as $\sum^{+\infty}_{j=e_\cL r_i}a^{(i)}_j\rho^j$ for 
some sequence of $a^{(i)}_j\!\in\!A_\cL$ \cite[corollary, pg.\ 68]
{koblitz}. Note in particular that $\frac{x_i}{a^{(i)}\rho^{e_\cL r_i}}$ 
thus lies in $R_p\!\setminus\!M_p$ for any 
$a^{(i)}\!\in\!A_\cL\!\setminus\!M_p$. 

Theorem \ref{thm:lip} thus implies that the number of isolated roots 
$(x_1,\ldots,x_n)$ of $F$ in $(\C^*_p)^n$ satisfying 
$(\ord_p x_1,\ldots,\ord_p x_n)\!=\!r$ {\bf and} 
$\frac{x_1}{a^{(1)}\rho^{e_\cL r_1}}\!\equiv \cdots \equiv\!\frac{x_n}{
a^{(n)}\rho^{e_\cL r_n}}\!\equiv\!1 \ (\mod \ M_p)$ is no more than 
$C_p(m,n,(1/e_\cL,\ldots,1/e_\cL))$). Furthermore, since 
$M_p\cap \cL\!\subset\!M_p$, 
we obtain the same statement if we restrict to roots in $(\cL^*)^n$ and use 
congruence $\mod \ M_p\cap \cL$ instead. 

Since there are $q_\cL-1$ possibilities for each $a^{(i)}_0$,  
our last observation tells us that the number of isolated roots 
$(x_1,\ldots,x_n)$ of $F$ in $(\cL^*)^n$ satisfying 
$(\ord_p x_1,\ldots,\ord_p x_n)\!=\!r$ is no more than 
$(q_\cL-1)^nC_p(m,n,(1/e_\cL,\ldots,1/e_\cL))$. So the total number of 
isolated roots of $F$ in $(\cL^*)^n$ is no more than 
$u(m,n)(q_\cL-1)^n C_p(m,n,(1/e_\cL,\ldots,1/e_\cL))$. Since $e_\cL\!\leq\!d$ 
and $q_\cL\!\leq\!p^d$, an elementary calculation yields our desired bound. 
\qed 

A simple consequence of our last proof is that 
there is a natural injection from the set of possible valuation 
vectors of an isolated root of $F$ to the set of lower facets\footnote{ 
cf.\ part 3 of remark \ref{rems:convex} of section \ref{sec:id}. Recall 
that a facet of a $d$-dimensional polytope is simply a face of dimension 
$d-1$. } of a particular  polytope. In particular, we can define 
$\pmb{\widehat{\Sigma}_p(F)}$ to be $\newt_p\left(\sum^k_{i=1}f_i\right)$
or the Minkowski sum $\sum^n_{i=1}\newt_p(f_i)$, according as $k\!>\!n$ or 
$k\!=\!n$, and immediately obtain the following corollary. 
\begin{cor} 
\label{cor:local}
Following the notation above, we have 
\[B(\cL,m,n)\leq\cF(F)(q_\cL-1)^nC_p(m,n,(1/e_\cL,\ldots,1/e_\cL)),\] 
where $\cF(F)$ is the number of lower facets 
of $\widehat{\Sigma}_p(F)$, and $q_\cL$ and $e_\cL$ are respectively 
the residue field cardinality and ramification index of $\cL$. \qed 
\end{cor}  

\begin{ex}
\label{ex:3m2}
Returning to example \ref{ex:3m}, observe that $f_1$ has 
$\leq\!3$ monomial terms (so $\newt_2(f_1)$ has $\leq\!3$ edges) and 
$f_2$ has $\leq\!\mu$ monomial terms (so $\newt_2(f_2)$ has 
$\leq\!2\mu$ edges (cf.\ our use of Euler's formula in the proof of the 
local case of theorem \ref{thm:big})). So we in fact have 
$\cF(F)\!\leq\!6\mu$ (for all $\mu\!\geq\!4$) and $\cF(F)\!\leq\!9$ 
(for $\mu\!=\!3$). Corollary \ref{cor:local} then implies improved upper 
bounds of 
\[ 304(\mu-1)\mu\left(1+\log_2\left(\frac{\mu-1}{0.693}\right)\right) 
\text{ (for all } \mu\!\geq\!4 \text{) \ \ \ \ \ and \ \ \ \ \ } 2304 
\text{ (for } \mu\!=\!3 \text{)} \] 
for the number of roots of $F$ in $(\Q^*_2)^2$. 
Also, our refined bound is smaller than the aforementioned real analytic bound 
(cf.\ example \ref{ex:3m} of section \ref{sec:intro}) for all $m\!\geq\!17$, 
where the two bounds begin to exceed $456800$. \dia 
\end{ex}  
\begin{ex} 
It is entirely possible that the maximum number of roots in 
$(\cL^*)^n$ of an $m$-sparse $n\times n$ polynomial system over 
$\cL$ is actually {\bf larger} for $\cL\!=\!\Q_2$ than for $\cL\!=\!\R$, 
for {\bf small} $m$ and $n$. In particular, a univariate trinomial over $\R$ 
clearly has at most $4$ nonzero real roots. 
However, $3x^{10}_1+x^2_1-4$ has exactly  
$6$ nonzero roots in $\Q_2$ and this is the maximum possible number 
of roots in $\Q^*_2$ for univariate trinomials over $\Q_2$ 
\cite[prop.\ 9.2]{lenstra2}. \dia 
\end{ex} 
 
\section{The Global Case of Theorem \ref{thm:big}} 
\label{sec:global} 
Let us start with a construction from \cite[sec.\ 8]{lenstra2}  
for the univariate case: First, fix a group homomorphism $\Q \longrightarrow 
\C^*_2$, written $r\mapsto 2^r$, with the property that $2^1\!=\!2$. 
To construct $2^r$ for an arbitrary rational $r$, choose $2^{1/n!}$ 
inductively to be an $n^\thth$ root of $2^{1/(n-1)!}$, and then define 
$2^{a/n!}$ to be the $a^\thth$ power of $2^{1/n!}$ 
for any $a\!\in\!\Z$. Clearly, $\ord_2(2^r)\!=\!r$ for each 
$r\!\in\!\Q$. For $j,e\!\in\!\N$ we then define the subgroups 
$U_e$ and $T_j$ of $\C^*_2$ by $U_e\!:=\!\{x \; | \; \ord_p(x-1)\!\geq\!1/e\}$ 
and $T_j\!:=\!\{\zeta \; | \; \zeta^{2^j-1}\!=\!1\}$. Note that  
$U_e\subseteq U_{e'}$ if $e\!\leq\!e'$, and $T_j\!\subseteq\!T_{j'}$ 
if $j$ divides $j'$. 

What we now show is that in addition to having few roots of bounded 
degree over $\Q_2$, $F$ has few roots in another suprisingly large piece of 
$(\C^*_2)^n$.  
\begin{lemma}
\label{lemma:len1} 
Let $e,j,k\!\in\!\N$, and let $F$ be an $m$-sparse $n\times n$ 
polynomial system over $\C_2$. Then $F$ has at 
most $\cF(F)(2^j-1)^nC_2(m,n,(1/e,\ldots,1/e))$ roots in the subgroup 
$(2^\Q\cdot T_j\cdot U_e)^n$ of $(\C^*_2)^n$, where $\cF(F)$ is 
as defined in corollary \ref{cor:local} of section \ref{sec:local}.  
\end{lemma} 

\noindent
{\bf Proof:}
First note that the case $n\!=\!1$, in slightly different notation, 
is exactly lemma 8.2 of \cite{lenstra2}. The proof there generalizes quite 
easily to our higher-dimensional setting. Nevertheless, for the convenience 
of the reader, let us give a succinct but complete proof. 

First note that by theorem \ref{thm:lip}, $F$ has no more than 
$C_2(m,n,(1/e,\ldots,1/e))$ roots in $U^n_e$. By the change 
of variables $(x_1,\ldots,x_n)\mapsto (\alpha_1y_1,\ldots,\alpha_ny_n)$ 
we then easily obtain the same upper bound for the number of roots 
of $F$ in any coset of $U^n_e$. Since $T^n_j$ clearly has order 
$(2^j-1)^n$, $F$ thus has no more than $(2^j-1)^n C_2(m,n,(1/e,\ldots,1/e))$ 
roots in any coset $(2^{r_1}T_jU_e)\times \cdots \times (2^{r_n} 
T_jU_e)$. Since Smirnov's Theorem implies, via our proof of the 
local case of theorem \ref{thm:big} (cf.\ section \ref{sec:local}), 
that a root $x\!\in\!(\C^*_2)^n$ of $F$ can produce no more than $\cF(F)$ 
possible distinct values for $(r_1,\ldots,r_n)\!:=\!(\ord_2 x_1, 
\ldots,\ord_2 x_n)$, we are done. \qed 

To at last prove the global case of theorem \ref{thm:big}, let us 
quote another useful result of Hendrik W.\ Lenstra, Jr.\; 
Recall that $\lceil x \rceil$ is the least integer greater than $x$. 
\begin{lemma}
\label{lemma:len2} 
\cite[lemma 8.3]{lenstra2}
Let $n\!\in\!\N$ and let $L$ be an extension of $\Q_2$ of 
degree $\leq\!D$. Then there is a $j\!\in\!\{1,\ldots,D\}$ such that 
$L^*\!\subseteq\!2^\Q T_j U_{\lceil d/j\rceil d}$. \qed 
\end{lemma} 

\noindent 
{\bf Proof of the Number Field Case of Theorem \ref{thm:big}:}\\ 
Since $\Q$ naturally embeds in $\Q_2$, we can assume 
$\cL$ is a subfield of $\C_2$ of finite degree over $\Q_2$. Then every root 
of $F$ in $(\C^*_2)^n$ of degree $\leq\!\delta$ over $\cL$ lies in 
$(L^{'*})^n$, where $L'$ is an extension of $\Q_2$ of degree at most 
$D\!:=\!d\delta$. So by lemma \ref{lemma:len2}, any such root of $F$ also 
lies in $\bigcup^D_{j=1} (2^\Q T_j U_{\lceil D/j\rceil D})$. 

From lemma \ref{lemma:len1} it now follows that the number of roots of $F$ of 
degree $\leq\!\delta$ over $\cL$ is no more than  
$\sum^{D}_{j=1}\cF(F)(2^j-1)^n
C_2\left(m,n,\left(\frac{1}{\lceil D/j\rceil D},\ldots,\frac{1}
{\lceil D/j\rceil D}\right) \right)$. 
Since $2^j-1\!\leq\!2^j$, $\cF(F)\!\leq\!u(m,n)$ 
(cf.\ the proof of the local case of theorem \ref{thm:big} 
in section \ref{sec:local}), and $C_2(m,n,(r,\ldots,r))$ is a 
decreasing function of $r$, we thus obtain by geometric 
series that $A(\cL,\delta,m,n)\leq 2^{nd\delta+1}u(m,n)C_2\left(m,n,
\left(\frac{1}{d^2\delta^2},\ldots,\frac{1}{d^2\delta^2}\right)\right)$. 
So by theorem \ref{thm:lip} and an elementary calculation we are done. \qed 

By leaving the last sum in our proof above unsimplified, we immediately obtain 
the following improvement of theorem \ref{thm:lip}. 
\begin{cor} 
\label{cor:global}
We have 
\[A(\cL,\delta,m,n)\leq \cF(F)\sum^{d\delta}_{j=1}
(2^j-1)^n C_2\left(m,n,\left(\frac{1}{\lceil d\delta/j\rceil d\delta},
\ldots,\frac{1}{\lceil d\delta/j\rceil d\delta}\right)\right), \] 
where $\cF(F)$ is as defined in corollary \ref{cor:local} of 
section \ref{sec:local}. \qed 
\end{cor}  

\section{Proving Theorem \ref{thm:lip}} 
\label{sec:lip} 

We begin with a clever observation of 
Hendrik W.\ Lenstra, Jr.\ on binomial coefficients, factorials, and least 
common multiples. Recall that $a|b$ means that $a$ and $b$ are 
integers with $a$ dividing $b$, and that $\lfloor x\rfloor$ denotes 
the greatest integer $\leq\!x$. 
\begin{dfn}
\label{dfn:dk} 
\cite[sec.\ 2]{lenstra2} 
For any nonnegative integers $m$ and $t$ define $\pmb{d_m(t)}$ 
to be the least common multiple of all integers that 
can be written as the product of at most $m$ pairwise 
distinct positive integers that are at most $t$ (and  
set $d_m(t)\!:=\!1$ if $m\!=\!0$ or $t\!=\!0$). 
Finally, for any $a\!\in\!\Z$, let us define $\pmb{\begin{pmatrix} a\\ 
t\end{pmatrix}}\!:=\!\prod^{t-1}_{i=0}\frac{a-i}{t-i}$ (and 
set $\begin{pmatrix}a\\ 0\end{pmatrix}\!:=\!1$). \dia 
\end{dfn} 
\begin{lemma} 
\cite[sec.\ 2]{lenstra2} 
\label{lemma:lenstra} 
Following the notation of definition \ref{dfn:dk}, we have...  
\begin{itemize}
\item[{\bf (a)}]{ $d_m(t)|n!$}
\item[{\bf (b)}]{ $m\!\geq\!t \Longrightarrow d_m(t)\!=\!t!$} 
\item[{\bf (c)}]{ $0\!\leq\!i\!\leq\!m\!<\!t 
\Longrightarrow i!|d_m(t)$}
\item[{\bf (d)}]{ $t\!\geq\!1 \Longrightarrow 
\ord_p d_m(t)\!\leq\!m\lfloor \log_p t\rfloor$} 
\end{itemize} 
\mbox{Furthermore, if $A\!\subset\!\Z$ 
is any set of cardinality $m$, then there are rational numbers}\\ 
\mbox{$\gamma_0(A,t),\ldots,\gamma_{m-1}(A,t)$ such that:}  
\begin{enumerate}
\item{the denominator of $\gamma_j(A,t)$ divides $d_{m-1}(t)/j!$ if 
$t\!\geq\!m$ and $\gamma_j(A,t)\!=\!\delta_{jt}$ 
otherwise.\footnote{$\delta_{ij}$ denoting the 
{\bf Kronecker delta}, which is $0$ when $i\!\neq\!j$ and $1$ when 
$i\!=\!j$.} }
\item{$\begin{pmatrix}a\\ t\end{pmatrix}=\sum^{m-1}_{j=0} 
\gamma_j(A,t)\begin{pmatrix}a\\ j\end{pmatrix}$ for all $a\!\in\!A$. \qed} 
\end{enumerate}
\end{lemma} 

Our proof of theorem \ref{thm:lip} will consist of a careful application  
of Smirnov's Theorem to the ``shifted'' polynomial system 
$G(x_1,\ldots,x_n)\!:=\!F(1+x_1,\ldots,1+x_n)$. (So roots of 
$F$ close to $(1,\ldots,1)$ are simply translations of roots of $G$ 
close to $(0,\ldots,0)$.) Since the $g_i$ can be highly non-sparse, one 
might not expect Smirnov's Theorem to give bounds independent of 
the degrees of the $f_i$ on the number of roots of $G$ close $(0,\ldots,0)$. 
However, lemma \ref{lemma:lenstra}, and lemmata 
\ref{lemma:ineq} and \ref{lemma:mink} below, save the day.  

\begin{lemma}
\label{lemma:ineq} 
Let $\pmb{c}\!:=\!\frac{e}{e-1}$ (so $c\!\leq\!1.582$) and 
$t_1,r_1,\ldots,t_n,r_n\!>\!0$. Then\\ 
\scalebox{.93}[1]{$\sum\limits^n_{i=1}(r_it_i-(m-1)\log_p t_i)\leq (m-1)
\sum\limits^n_{i=1}r_i \Longrightarrow \sum\limits^n_{i=1}
r_it_i\leq c(m-1)\left[\left(\sum\limits^n_{i=1} r_i\right) 
+ \log_p\left(\frac{(m-1)^n}{r_1\cdots r_n\log^n p}\right)\right]$.}  
\end{lemma} 

\noindent
{\bf Proof:} Here we make multivariate extensions of some observations of 
Lenstra from \cite[prop.\ 7.1]{lenstra2}: 
First note that it is easily shown via basic calculus that 
$1-\frac{\log x}{x}$ assumes its minimum (over the positive 
reals), $1/c$, at $x\!=\!e$. So for all $x\!>\!0$ we have $x\!\geq\!(\log x) 
+x/c$. Letting $t,r\!>\!0$, $w\!:=\!\frac{m-1}{r\log p}$, and 
$x\!:=\!t/w$, we then obtain\\ 
\scalebox{.85}[1]{
$rt\!\geq\!rwx\!\geq\!rw((\log x)+x/c)\!=\!rw(\log t)
-rw(\log w)+rt/c\!=\!(m-1)(\log_p t)-(m-1)\log_p\left(\frac{m-1}{r\log p}
\right)+rt/c$. 
} 
Substituting $r\!=\!r_i$, $t\!=\!t_i$, and summing over $i$ then implies 
\[ (\star) \ \ \ \  \sum\limits^n_{i=1} r_it_i 
\geq (m-1)\left(\sum\limits^n_{i=1} 
\log_p t_i\right)-(m-1)\log_p
\left(
\frac{(m-1)^n}{r_1\cdots r_n\log^n p}
\right)
+\frac{1}{c}\sum\limits^n_{i=1}r_it_i. \] 

Now suppose that 
\[ (\star\star) \ \ \ \ 
\sum\limits^n_{i=1}
r_it_i> c(m-1)\left[\left(\sum\limits^n_{i=1} r_i\right) 
+ \log_p\left(\frac{(m-1)^n}{r_1\cdots r_n\log^n p}\right)\right]. \] 

Substituting ($\star\star$) into the {\bf last} sum of the 
{\bf right} hand side of our inequality ($\star$) then tells us\\ 
\scalebox{.9}[1]{$\sum\limits^n_{i=1} r_it_i\!>\! 
(m-1)\left(\sum\limits^n_{i=1} 
\log_p t_i\right)-(m-1)\log_p\left(\frac{(m-1)^n}
{r_1\cdots r_n\log^n p}\right)+ 
(m-1)\left[\left(\sum\limits^n_{i=1} r_i\right) 
+ \log_p\left(\frac{(m-1)^n}{r_1\cdots r_n\log^n p}\right)\right]$.}\\ 
So we obtain $\sum\limits^n_{i=1} r_it_i 
> (m-1)\left(\sum\limits^n_{i=1} 
\log_p t_i\right) +(m-1)\left(\sum\limits^n_{i=1} r_i\right)$, which 
can be rearranged into   
\[ (\star\star\star) \ \ \ \ \sum\limits^n_{i=1} (r_it_i-(m-1)\log_p t_i) 
> (m-1)\sum\limits^n_{i=1}r_i. \] 
So ($\star\star$) $\Longrightarrow$ ($\star\star\star$), 
and we conclude simply by taking the contrapositive. \qed 

The following lemma is a simple consequence of the basic 
properties of polytopes, their faces, and their mixed volumes \cite{buza}. 
\begin{lemma}
\label{lemma:mink}
Following the notation of \ref{sec:id}, let $G\!:=\!(g_1,\ldots,g_n)$ be any 
$n\times n$ polynomial system and let $r\!:=\!(r_1,\ldots,r_n)$ be 
such that $r_i\!>\!0$ for all $i$. Also let 
$\pmb{w(g_i,r)}\!:=\pi\left( \bigcup\limits_{\substack{
\hat{s}:=(s_1,\ldots,s_n,1)\\ 
s_i\geq r_i \text{ for all } i}} \! \! \! \! \! \! \! \! \! 
\newt^{\hat{s}}_p(g_i)\right)$ for all $i$. Then  
$\sum\limits_{\substack{
\hat{s}:=(s_1,\ldots,s_n,1)\\ 
s_i\geq r_i \text{ for all } i}} 
\cM(\pi(\newt^{\hat{s}}_p(G)))\!\leq\!\cM(\conv(w(g_1,r)),\ldots,
\conv(w(g_n,r)))$. In particular, if $Q\!=\!\{(t_1,\ldots,t_n)\!\in\!\Rn \; 
| \; r_1t_1+\cdots r_nt_n\!\leq\!1 \text{ and } t_j\!\geq\!0 \text{ 
for all } j\}$, then 
$\cM(\underset{n \text{ times}}{\underbrace{Q,\ldots,Q}})\!=\!1
\left/\prod^n_{i=1}r_i\right.$. \qed 
\end{lemma} 

\noindent 
{\bf Proof of Theorem \ref{thm:lip}:}\mbox{}\\ 
First note that just as in the proof of the local case of theorem 
\ref{thm:big} (cf.\ section \ref{sec:local}), we have that 
$k\!<\!n$ or $m\!\leq\!n$ implies that there are no isolated roots 
whatsoever. So we can assume that $m\!\geq\!n+1$. Also, 
again like in the proof of the local case of theorem \ref{thm:big}, we can 
safely assume via lemma \ref{lemma:gh} that $k\!=\!n$ and observe that root 
multiplicities are preserved if we already had $k\!=\!n$ in our original input. 
Furthermore, if $m_i$ is the number of monomial terms occuring in $f_i$ 
for all $i$, then it is easily checked that $m_i\!\leq\!1$ for any $i$ implies 
that there are no isolated roots at all. So we can also assume that 
$m_1,\ldots,m_n\!\geq\!2$. 

Let us now set $g_i(x_1,\ldots,x_n)\!:=\!f_i(1+x_1,\ldots,1+x_n)$ for all $i$ 
and $G\!:=\!(g_1,\ldots,g_n)$. It is then clear that the number of 
isolated roots of $F$ with $\ord_p(x_i-1)\!\geq\!r_i$ for all $i$ is the 
same as the number of isolated roots of $G$ with $\ord_p x_i \!\geq\!r_i$ 
for all $i$, and multiplicities are preserved by this change of variables. 
Smirnov's Theorem tells us that the latter number (counting multiplicities) 
is exactly\footnote{
Note that the sum over $s$ is actually infinite, but 
has only finitely many nonzero summands. This is because  
any polytope has only finitely many inner facet normals with last 
coordinate $1$, and it is only these terms which can possibly 
be nonzero.}  
$\sum\limits_{\substack{\hat{s}:=(s_1,\ldots,s_n,1)\\
s_i\geq r_i \text{ for all } i}} 
\cM(\pi(\newt^{\hat{s}}_p(G)))$.

Now let us define the following scaled standard simplex: 

\vspace{.5cm}
\noindent
\scalebox{.88}[1]
{$\pmb{S(m,n,r)}\!:=\!\left\{(t_1,\ldots,t_n)\!\in\!\Rn 
\; \left| \; \sum\limits^n_{j=1}
r_jt_j\leq c(m-1)\left[\left(\sum\limits^n_{j=1} r_j\right)
+ \log_p\left(\frac{(m-1)^n}{r_1\cdots r_n\log^n p}\right)\right]
\text{ and } t_j\!\geq\!0 \text{ for all } j\right.\right\}$.} 

\vspace{.5cm}
\noindent
Note then that by lemma \ref{lemma:mink}, 
$\cM(\underset{n \text{ times}}{\underbrace{S(m_1,n,r),\ldots,S(m_n,n,r)}})$ 
is exactly \[c^n\left. \prod^n_{i=1} 
\left\{(m_i-1)\left[\left(\sum^n_{j=1} r_j\right) + 
\log_p\left(\frac{(m_i-1)^n}{r_1\cdots r_n\log^n p}\right)\right]\right/ r_i
\right\},\] since mixed volume is multihomogeneous with respect to 
scalings \cite{buza}. Since $S(m,n,r)$ is clearly always convex, and since 
$w(g_i,r)$ is a union of convex hulls of subsets of $\supp(g_i)$, we also 
have that $w(g_i,r)\cap\supp(g_i)\!\subseteq\!S(m_i,n,r) \Longrightarrow 
\conv(w(g_i,r))\!\subseteq\!S(m_i,n,r)$ for all $i$. 

Since mixed volume is monotonic with respect to containment \cite{buza},  
lemma \ref{lemma:mink} then clearly implies that... 
\begin{center}
{\bf \scalebox{.95}[1]{To prove theorem 
\ref{thm:lip}, we need only show that $\pmb{w(g_i,r)\cap 
\supp(g_i)\!\subseteq\!S(m_i,n,r)}$ for all 
$\pmb{i}$.}} 
\end{center}
To do this, we will first prove that the valuations of the 
coefficients of any $g_i$ satisfy a ``slow decay'' condition, 
and then use convexity of the gently sloping lower faces of the 
$p$-adic Newton polytopes $\newt_p(g_i)$ to prove that $w(g_i,r)\cap 
\supp(g_i)\!\subseteq\!S(m_i,n,r)$ for all $i$. 

Let us temporarily abuse notation slightly to avoid a profusion of indices 
and respectively write $f$, $g$, and $m$ in place of $f_i$, $g_i$, and 
$m_i$ (for some arbitrary fixed $i$). Letting $D_i\!:=\!\deg_{x_i} f$, it is 
clear that we can write $g(x)\!:=\!\sum_{j\in \prod^n_{i=1}\{0,\ldots,D_i\}} 
b_jx^j$, where $\pmb{b_j}\!:=\!\sum_{a\in A} c_a \prod^n_{i=1}
\begin{pmatrix}a_i\\ j_i\end{pmatrix}$, 
$f(x)\!=\!\sum\limits_{a=(a_1,\ldots,a_n)\in A}\pmb{c_a} x^a$ (with every 
$c_a$ nonzero), $j\!=\!(j_1,\ldots,j_n)$, and $\pmb{A}\!:=\!\supp(f)$. 
Since $f\!\neq\!0$ we have $g\!\neq\!0$ and thus not all the $b_j$ vanish. 
Note also that $D_i\!=\!0 \Longrightarrow \supp(g)\!\subseteq\!\{ 
x\in\!\Rn \; | \; x_i\!=\!0\}$. Letting $\pi_i : \Rn \longrightarrow \R^{n-1}$ 
denote the natural orthogonal projection forgetting the $i^\thth$ 
coordinate, it is then clear that 
\[\pi_i\left(S(m,n,(r_1,\ldots,r_n))\!\cap\!\{x\!\in\!\Rn \; 
| \; x_i\!=\!0\}\right)\!\supseteq\!S(m,n-1,
(r_1,\ldots,r_{i-1},r_{i+1},\ldots,r_n)),\]  
and thus $D_i\!=\!0$ implies that we can reduce to a case where $n$ is 
smaller. So we can assume henceforth that $D_1,\ldots,D_n\!>\!0$. 

By lemma \ref{lemma:lenstra} there are rational numbers 
$\{\gamma^{(i)}_j(t_i)\}$, with 
\mbox{$(i,j)\in\{1,\ldots,n\}\times\{0,\ldots,m-1\}$,}  
such that for all $a\!=\!(a_1,\ldots,a_n)\!\in\!A$ we have 
$\begin{pmatrix}a_i\\ t_i\end{pmatrix}\!=\!\sum^{m-1}_{j=0}
\gamma^{(i)}_j(t_i) \begin{pmatrix}a_i\\ j\end{pmatrix}$  
and the denominators of the $\{\gamma^{(i)}_j(t_i)\}$ are not too 
divisible by $p$. (We will make the latter assertion precise in a moment.) 

We thus obtain that for all $t\!:=\!(t_1,\ldots,t_n)\!\in\!\prod\limits^n_{i=1} 
\{0,\dots,D_i\}$,  
\[ b_t=\sum_{a\in A} c_a 
\prod^n_{i=1} \begin{pmatrix}a_i\\ t_i\end{pmatrix} 
 = \sum_{a\in A} c_a \prod^n_{i=1} 
\sum^{m-1}_{j_i=0} \left(\gamma^{(i)}_{j_i}(t_i)\begin{pmatrix}a_i\\ 
j_i\end{pmatrix}\right) 
=\sum_{a\in A} c_a \!\!\!\!
\sum_{j\in\{0,\ldots,m-1\}^n} \prod^n_{i=1}\left( 
\gamma^{(i)}_{j_i}(t_i)\begin{pmatrix}a_i\\ 
j_i\end{pmatrix}\right) \]
\[ = \!\!\!\!\! \sum_{j\in\{0,\ldots,m-1\}^n} 
\left(\prod^n_{i=1}\gamma^{(i)}_{j_i}(t_i) \right) 
\sum_{a\in A} c_a \prod^n_{i=1}\begin{pmatrix}a_i\\ 
j_i\end{pmatrix} 
= \!\!\!\!\! \sum_{j\in\{0,\ldots,m-1\}^n} 
\left(\prod^n_{i=1}\gamma^{(i)}_{j_i}(t_i)\right) b_j. \] 
So the coefficients $\{b_t\}_{t\in\prod^n_{i=1}\{0,\ldots,D_i\}}$ of $g$ are 
completely determined by a {\bf smaller} set of coefficients corresponding to 
the exponents of $g$ lying in $\{0,\ldots,m-1\}^n$. Even better, 
lemma \ref{lemma:lenstra} tells us that $t_i\!\leq\!m-1 \Longrightarrow 
\gamma^{(i)}_{j_i}(t_i)\!=\!0$ for all $j_i\!\neq\!t_i$. So we in fact 
have that 
\[ (\heartsuit) \ \ \ \ t_i\!\leq\!m-1 \Longrightarrow \text{ the 
recursive sum for } b_t \text{ has {\bf no} terms corresponding to any } 
j \text{ with } j_i\!\neq\!t_i. \] 

Given this refined recursion for $b_t$ we can then derive that 
the $p$-adic valuation of $b_t$ decrease slowly and in a highly 
controlled manner: First note 
that our recursion, combined with ($\heartsuit$) and the ultrametric 
inequality, implies that 
\[ (\star) \ \ \ \ \ \ord_p b_t \geq\!\min_{j\in M_t} 
\left\{ \ord_p(b_j)+\sum^n_{i=1}\ord_p\gamma^{(i)}_{j_i}(t_i)
\right\} \text{ for all } t\!\in\!\prod^n_{i=1}\{0,\ldots,D_i\},\]
where $M_t$ is the subset of $\{0,\ldots, m-1\}^n$ obtained by 
the intersection, over all $i$ with \mbox{$t_i\!\leq\!m-1$,} of the 
hyperplanes $\{t\!\in\!\Rn \; | \; t_i\!=\!j_i\}$. 
Then, by the definition of a face with inner normal $(s,1)$, 
we have $(t,b_t)\!\in\!\newt^{(s,1)}_p(g) 
\Longrightarrow \left(\sum^n_{i=1}s_it_i\right)+\ord_pb_t\!\leq\! 
\left(\sum^n_{i=1}s_ij_i\right)+ \ord_pb_j$
for all $j\!\in\!\prod^n_{i=1}\{0,\ldots,D_i\}$. So 
for all such $j$ we must have $\ord_p b_j \geq \ord_p b_t + 
\sum^n_{i=1}s_i(t_i-j_i)$. In particular, we obtain that   
\[ (\star\star) \ 
[(t,b_t)\in\newt^{(s,1)}_p(g) \text{ and } 
t_i\!\geq\!j_i \text{ and } s_i\!\geq\!r_i \text{ for all } i] 
\Longrightarrow \ord_p b_j \geq \ord_p b_t + \sum^n_{i=1}r_i(t_i-j_i).\] 

Since $t\!\in\!\supp(g)$ and $(t,\ord_p b_t)\!\in\!\newt^{(s,1)}_p(g)$  
implies that $\ord_p b_t\!<\!\infty$, we can thus combine ($\star$) and 
($\star\star$) to obtain that \[t\!\in\!w(g,r)\cap\supp(g) \Longrightarrow 
\ord_p b_t\!\geq\!\min\limits_{j\in M_t}
\left\{ \ord_p (b_t) + \sum^n_{i=1}\left(r_i(t_i-j_i)+
\ord_p\gamma^{(i)}_{j_i}(t_i)\right)\right\}.\] 

Cancelling and rearranging terms, we thus obtain that 
\mbox{$t\!\in\!w(g,r)\cap\supp(g) \Longrightarrow$}
\[ \sum^n_{i=1}r_it_i 
\leq \max\limits_{j\in M_t}\left\{ \sum^n_{i=1} 
\left(j_ir_i -\ord_p(\gamma^{(i)}_{j_i}(t_i)
\right)\right\}\leq\max\limits_{j\in \{0,\ldots,m-1\}^n}\left\{ \sum^n_{i=1}
\left(j_ir_i -\ord_p(\gamma^{(i)}_{j_i}(t_i)
\right)\right\}. \] 
Since lemma \ref{lemma:lenstra} tells us that 
$-\ord_p\gamma^{(i)}_{j_i}(t_i)\!\leq\!(m-1)(\log_p 
t_i)-\ord_p(j_i!)$ for all $i$, we then obtain\\ 
\scalebox{.95}[1]{$(\clubsuit) \ \ \ \ \  \sum^n_{i=1}\left(r_it_i - 
(m-1)\log_p t_i\right) \leq \max\limits_{j\in\{0,\ldots,m-1\}^n}\left\{ 
\sum^n_{i=1} \left(j_ir_i -\ord_p(j_i!) \right)\right\}\leq(m-1)
\sum^n_{i=1}r_i$.}\\ 
So by lemma \ref{lemma:ineq} we obtain that 
$w(g,r)\cap\supp(g)\!\subseteq\!S(m,n,r)$, and thus 
$w(g_i,r)\cap\supp(g_i)\!\subseteq\!S(m_i,n,r)$ for all $i$. \qed 

\section{Connections to Complexity Theory} 
\label{sec:pnp}
Thanks to our results, we now know in particular that the maximum 
number of isolated rational roots of a $k\times n$ polynomial system over 
$\Q$ depends polynomially on the number of distinct exponent vectors, for 
fixed $n$. Here we note that it would be of considerable interest to know if 
this polynomiality persists relative to even more efficient encodings of 
polynomials. 

In particular, instead of monomial expansions (a.k.a.\ the 
{\bf sparse encoding}), consider the {\bf 
straight-line program (SLP) encoding} for 
a univariate polynomial \cite[sec.\ 7.1]{bcss}: That is, suppose we have 
$p\!\in\!\Z[x_1]$ expressed as a sequence of the form 
$(1,x_1,q_2,\ldots,q_N)$, where $q_N\!=\!p$ and for all $i\!\geq\!2$ we 
have that $q_i$ is a sum, difference, or product of some pair of elements 
$(q_j,q_k)$ with $j,k\!<\!i$. 
Let $\pmb{\tau(p)}$ denote the smallest possible value of $N-1$, i.e., 
the smallest length, for such a computation of $p$. Clearly, 
$\tau(p)$ is no more than the number of monomial terms of $p$, and 
is often dramatically smaller. 
\begin{thm}
\label{thm:smale} 
\cite[thm.\ 3, pg.\ 127]{bcss} 
Suppose there is an absolute constant $\kappa$ such that 
for all nonzero $p\!\in\!\Z[x_1]$, the 
number of distinct roots of $p$ in $\Z$ is no more than $(\tau(p)+1)^\kappa$. 
Then $\pp_\C\!\neq\!\np_\C$. \qed  
\end{thm} 
In other words, an analogue (regarding complexity 
theory over $\C$) of the famous unsolved $\pp\stackrel{?}{=}\np$ question from 
computer science (regarding complexity theory over the ring $\Z/2\Z$) would 
be settled. The question of whether $\pp_\C\stackrel{?}{=}\np_\C$ remains 
open as well but it is known that $\pp_\C\!=\!\np_\C \Longrightarrow 
\np\!\subseteq\!\bpp$. (This observation is due to Steve Smale and 
was first published in \cite{shub}.) The complexity class $\bpp$ is central 
in randomized complexity and the last inclusion (while widely disbelieved) 
is also an open question. The truth of the hypothesis of theorem 
\ref{thm:smale}, also know as the {\bf $\pmb{\tau}$-conjecture}, is 
yet another open problem, even for $\kappa\!=\!1$. 

One can reasonably suspect that a sufficiently good upper bound for the number 
of integral roots of an $m$-sparse $k\times n$ polynomial system could be 
applied to settling the $\tau$-conjecture:\footnote{ 
i.e., Diophantine results for multivariate polynomial systems in the 
sparse encoding can be useful for Diophantine problems involving univariate 
polynomials in the SLP encoding.}  Indeed, 
any computation of $p$ of length $\tau(p)$ can be specialized to obtain a 
computation of the same length for $p(k)$ for any integer $k$. Finding an 
integral root of $p$ can then be reinterpreted as finding all possible values 
for the first coordinate of an integral root of the 
$2\tau(p)$-sparse $\tau(p)\times 
(\tau(p)+1)$ polynomial system defined by the corresponding computational 
sequence for $p$. However, the number of variables grows linearly with 
$\tau(p)$, so this route toward an application of theorem \ref{thm:big} would 
at best give us an upper bound exponential in $\tau(p)$. 
For better or worse, we thus arrive at a Diophantine problem currently out of 
our grasp: finding sharp bounds on the number of integral points on certain 
algebraic sets defined by quadratic binomials and linear trinomials. 

A reasonable alternative approach would be to use the embedding 
of $\Q$ in another complete field --- $\R$, in particular. 
Over $\R$ there are results for univariate polynomials involving an 
even sharper encoding: For any $p\!\in\!\R[x]$, let its {\bf additive 
complexity}, $\pmb{\sigma(p)}$, be the minimal number of additions and 
subtractions necessary to express $p$ as an elementary algebraic expression 
with {\bf constant} exponents. e.g., 
$p(x)\!=\!(1-(x+2)^{100})^{97}+243(x-7)^{999}$ 
has $\sigma(p)\leq\!4$, and it is clear that 
$\tau(p)\!\geq\!4$. More generally, it is easily checked that 
$\sigma(p)\!\leq\!\tau(p)$ for all $p\!\in\!\Z[x_1]$. Remarkably, one can 
bound the number of {\bf real} roots of $p$ solely in terms of $\sigma(p)$. 
This was known since the work of Allan Borodin and Stephen A.\ Cook around 
1974 \cite{bocook}, and the best current upper bound is Jean-Jacques Risler's 
$C^{\sigma(p)^2}$, for some absolute constant $C\!\in\!(1,32)$ 
\cite{grigo,risler}. Unfortunately, there are examples of $p\!\in\!\Z[x_1]$ 
with $\sigma(p)\!=\!\cO(r)$ and at least $2^r$ real roots (all of which 
are irrational) \cite[sec.\ 3, pg.\ 13]{four}. So additive complexity 
is too efficient an encoding to be useful in settling the $\tau$-conjecture, 
at least over $\R$. 

Whether analogous (hopefully polynomial) bounds in terms of a sharper 
encoding exist in our {\bf arithmetic} setting is an open question, even for 
$n\!=\!1$. In particular, it is interesting to note that the only obstructions 
to refining theorem \ref{thm:big} to a sharper encoding 
are (a) the strong dependence of the quantity $\cF(F)$, arising from our 
application of Smirnov's Theorem, on the number of monomial terms, and (b) the 
existence of an analogue of theorem \ref{thm:lip} for a sharper encoding. 
 
As for actually {\bf finding} all the isolated rational roots of $F$, 
there is both good news and bad news: The bad news is that 
one can {\bf not} have a polynomial time algorithm (relative to 
the sparse encoding) for $n\!>\!1$. The good news is that there 
{\bf is} a polynomial time algorithm (relative to the sparse encoding) 
for $n\!=\!1$, and that the counter-examples for $n\!>\!1$ are very simple. 

In particular, if we take $\cL\!=\!\Q$ and measure the input size simply 
as the number of digits needed to write the coefficients {\bf and} exponents 
of $F$ in, say, binary; then it possible for an isolated rational root of $F$ 
to have bit size\footnote{ The bit size of an integer is thus implicitly the 
number of digits in its binary expansion, and the bit size of a rational 
number can be taken as the maximum of the bit sizes of its numerator and 
denominator (written in lowest terms).} exponential in the bit size of $F$:
Simply consider $k\!=\!n\!=\!2$, $m\!=\!4$, and 
$F\!:=\!(x_1-x^D_2,x_2-2)$. This particular example clearly has bit size 
$\cO(\log D)$ but its one rational root $(2^D,2)$ has a first coordinate 
of bit size $D$ --- exponential in the bit size of $F$. Thus one can't 
even write the output in polynomial time relative to the sparse encoding. 
Similar examples with bit size $\cO(n\log D)$ and having a single rational 
root, but with root coordinates of bit size $\Omega(D^n)$, are easy to 
construct for all $n\!\geq\!3$ via the same recursive idea \cite[pg.\ 16, 
complication {\bf Q}$\pmb{_2}$]{four}. With a bit more 
work one can even show that such roots of ``excessively large'' bit size 
occur not only in a worst case sense but also in an {\bf average case} sense.

On the other hand, it is a fortunate accident that the {\bf absolute 
logarithmic height} of a complex root of $F$ of degree $\leq\!\delta$ over 
$\cL$ (and thus equivalently, the bit size of such a root) is polynomial in the 
bit size of $F$ for $n\!=\!1$ and $\cL$ a number field \cite[prop.\ 
2.3]{lenstra1}. This is what permits a clever polynomial time algorithm for 
solving $F$ when $n\!=\!1$ and $\cL$ and $\delta$ are fixed 
\cite[first theorem]{lenstra1}.\footnote{Lenstra's algorithm 
has complexity exponential in $\delta$ and the degree of $\cL$ over $\Q$. 
Note however that his algorithm is considerably faster than the 
well-known Lenstra-Lenstra-Lovasz factoring algorithm \cite{lll}: the 
latter algorithm would only solve $x^D+ax+b\!=\!0$ over the rationals 
in time {\bf exponential} in $\log D$. }  
For $n\!>\!1$ it thus appears that the only way to achieve a 
polynomial time algorithm would be to allow a more efficient encoding 
of the output than expanding into digits. 
In particular, it is an open question, even for $n\!=\!2$, whether one can 
always find {\bf SLP's}, of length polynomial in the bit size of $F$, for the 
isolated rational roots of $F$. 

Alternatively, one can simplify the question of solving and simply 
ask how many isolated rational roots $F$ has, or whether $F$ has any 
isolated rational roots at all. This is addressed in \cite[thms.\ 1.3 and 
1.4]{jcs}, where it is shown that the truth of the Generalized Riemann 
Hypothesis implies that detecting a strong form of {\bf non}-solvability over 
the rationals (transitivity of the underlying Galois group) can be done within 
the complexity class $\pp^{\np^\np}$, provided the underlying complex zero set 
is finite. In the latter result, $n$ is allowed to be part of the input 
and can thus vary. 

\section*{Acknowledgements} 
The author thanks Raphael Hauser and Gregorio Malajovich for useful 
discussions. 

\footnotesize
\bibliographystyle{acm}

\end{document}